\documentclass[letterpaper]{article}

\usepackage{aaai19}
\usepackage{times}
\usepackage{helvet}
\usepackage{natbib}
\usepackage{courier}
\usepackage{url}  
\usepackage{graphicx}  
\usepackage{amsfonts}       

\usepackage{amsmath,cases,amssymb,verbatim,amsthm}
\usepackage{algorithm2e,algorithmic}
\usepackage{mathtools}

\DeclareMathOperator*{\minimize}{minimize}
\DeclareMathOperator*{\maximize}{maximize}

\newcommand{\norm}[1] {\|#1\|}
\newcommand{\mdot}[2]{\langle#1,#2\rangle}

\def\bR{\mathbb{R}}

\AlgoDontDisplayBlockMarkers
\SetAlgoNoEnd
\SetAlgoNoLine
\title{Low-rank semidefinite programming for the MAX2SAT problem}

\author{
   Po-Wei Wang\\
   Machine Learning Department\\
   Carnegie Mellon University\\
   Pittsburgh, PA 15213 \\
   \texttt{poweiw@cs.cmu.edu}
   \And
   J. Zico Kolter\\
   School of Computer Science\\
   Carnegie Mellon University\\
   Pittsburgh, PA 15213, and \\
   Bosch Center for Artificial Intelligence\\
   Pittsburgh, PA 15222 \\
   \texttt{zkolter@cs.cmu.edu}
}

\begin{document}
\maketitle              

\setcounter{secnumdepth}{2}
\begin{abstract}
        This paper proposes a new algorithm for solving MAX2SAT problems based on combining search methods with semidefinite programming approaches.  Semidefinite programming techniques are well-known as a theoretical tool for approximating maximum satisfiability problems, but their application has traditionally been very limited by their speed and randomized nature.  Our approach overcomes this difficult by using a recent approach to low-rank semidefinite programming, specialized to work in an incremental fashion suitable for use in an exact search algorithm.  The method can be used both within complete or incomplete solver, and we demonstrate on a variety of problems from recent competitions. Our experiments show that the approach is faster (sometimes by orders of magnitude) than existing state-of-the-art complete and incomplete solvers, representing a substantial advance in search methods specialized for MAX2SAT problems.
\end{abstract}
\section{Introduction}

The maximum satisfiability (MAXSAT) task --- the optimization-based version of satisfiability, where the goal is to find a maximal set of clauses to be satisfied --- is a foundational problem in combinatorial optimization.  These problems can be written as the optimization task
\begin{equation}
\maximize_{v_1,\ldots,v_n\in \{0,1\}} \sum_{j=1}^m \bigvee_{i \in S^+_j} v_i \bigvee_{i \in S^-_j} \neg v_i,
\end{equation}
where $v_1,\ldots,v_n$ denote the binary optimization variables, and $S^{+}_j$/$S^{-}_j$ denotes the set of variables and negated variables respectively in the $j$th clause.  MAXSAT problems have been used to solve probabilistic inference \citep{park2002using}, planning \citep{zhang2012maxsat}, and clustering \citep{berg2017cost}. And the recent MAXSAT contests \citep{maxsatrace2016} have aimed at evaluating a wide number of different solution methods on a wide variety of both real and synthetic problems.

In this paper, we present an approach based upon low-rank semidefinite programming (SDP) and a simple branch-and-bound strategy for solving the simplest case of the MAXSAT problem: the MAX2SAT problem, where each clause contains only two variables.  While this is naturally a very limited special case of the MAXSAT problem, unlike the SAT setting (where 2-SAT can be solved in polynomial time), the MAX2SAT problem is NP-COMPLETE and forms a non-trivial starting point for more general MAXSAT problems.  Although SDP relaxations of the MAX2SAT problem have been known for some time \citep{goemans1994new,goemans1995improved,gomes2006power}, they have not generally been viewed as practical strategies for solving MAXSAT problems (of any type), owing to the high computational cost of solving SDPs via standard methods (e.g., interior point methods, whose runtime grows cubically in the number of variables).  Indeed, to the best of our knowledge, no entry in the history of the MAXSAT competition \citep{maxsatrace2016} has ever used an approach based upon semidefinite programming.
Furthermore, our approach also applies to more general MAXSAT problems, though the practical efficiency is much less for settings with a large number of variables per clause; but later in the paper, we will also highlight MAX3SAT results, for example, where the method is also still competitive.

In this work, we show that a recent method for low-rank semidefinite programming (the Mixing method \citep{wang2017mixing}) combined with branch-and-bound strategy, warm-starts, and simple pruning based upon a dual bound, can achieve state-of-the-art performance for MAX2SAT problems.  Specifically, for the MAXSAT 2016 competition problems\footnote{ The 2016 MAXSAT competition was the last year of the competition to features MAX2SAT instances, as the ``random'' problem track was removed in subsequent years.}, our single proposed approach dominates the best solvers (selecting the best solver for each individual problem instance for all the 2016 problems) in \emph{both} the incomplete and complete tracks of the competition.  While naturally these results need to be taken with a grain of salt, as the MAXSAT competition solvers of course did \emph{not} specialize for the MAX2SAT case, it is notable our single solver, using ultimately a very simple search procedure, can outperform highly specialized and tuned solvers so heavily in this domain across both competition tracks.  Thus, while a great deal of work remains to turn the approach into a general MAXSAT solver, the work here strongly suggests that semidefinite programming techniques can play a sizable role in the solution to MAXSAT problems going forward.

\section{Preliminaries and related work}

In this section we highlight a number of current approaches for solving MAXSAT problems, and review relevant recent work on semidefinite approximations and fast semidefinite programming methods.  We will specifically highlight the recent Mixing method \citep{wang2017mixing}, which forms the basis of our inner solution method.

\subsection{Discrete MAXSAT solvers}
The vast majority of modern methods for the MAXSAT problem are based upon discrete search method.  As our proposed approach in this work (which based upon continuous optimization) differs quite substantially from these approaches, we refer the reader to a recent survey \citep{morgado2013iterative} for much more detailed descriptions about the current state of the art.  However, broadly speaking, there have been two main classes for these discrete solvers: 1) those based upon bounding the solution via SAT method \citep{marques2011algorithms,koshimura2012qmaxsat,ansotegui2013sat,fu2006solving,le2010sat4j,een2006translating} which in turn exploit the heuristic developed by the SAT community such as those in the MiniSAT solver; these solvers typically are \emph{complete} in that they will both produce a satisfying assignment with some number of clauses satisfied \emph{and} a verification that this is the optimal solution to the problem.  And 2) those based upon local search \citep{luo2015ccls,luo2017ccehc}, which maintain and locally adjust a solution to satisfy an increasing number of clauses; these solvers typically are \emph{incomplete} in that they may quickly find an assignment, but often cannot prove whether or not it is optimal.  There are also aggregation-based solvers such as Open-WBO \citep{martins2014open}, which integrate several solvers for different instances.

\subsection{Approximation algorithms for MAXSAT}
In addition to exact solvers based upon combinatorial search, several approaches have also considered \emph{approximation} algorithms for the MAXSAT problem, both from theoretical and practical standpoints.  Generally speaking, these focus on obtaining a high approximation ratio $\alpha$, such that
\begin{equation}
        \text{Approximated objective} \geq \alpha \cdot \text{Optimal objective}.
\end{equation}
Specifically, \citet{goemans1994new} proposed an LP approximation algorithm for MAXSAT with $\alpha=0.75$.
They latter proposed also an SDP approximation \citep{goemans1995improved} for MAX2SAT with $\alpha=0.878$, which will form the basis for the SDP problem we solve here, and which we describe in more detail below.  \citet{feige1995approximating} used a different and more refined SDP to raise the approximation ratio to $\alpha=0.931$, and \citet{karloff19977} extend the result to deal with MAX3SAT with $\alpha=0.875$.
Interested readers can refer to \citet{biere2009handbook} for a comprehensive introduction.

Although several authors have noted the effectiveness of the SDP relaxation for the MAX2SAT problem \citep{lotgering2012sdp,gomes2006power} in terms of the approximation quality, no previous works we are aware of have succeeded at turning this into a practical algorithm.  Several efforts were made: \citet{anjos2004semidefinite} investigated the effectiveness of different SDPs for the MAXSAT problems, but didn't consider their integration into an exact solver. \citet{van2008sums} investigated the effectiveness of sum-of-squares based SDP solvers and found them not practical.
And \citet{rendl2010solving,krislock2017biqcrunch} demonstrated an SDP-based branch-and-bound framework for the related MAXCUT and binary quadratic problem (using interior point and bundle solvers) and created the Biq Church tool; however, this method was only feasible for small dense instances.

\subsection{Low-rank semidefinite programming}
As mentioned, the approach we ultimately use is one based upon semidefinite programming.
A general SDP can be written as the optimization problem
\begin{equation}
\begin{split}
\minimize_{X \in \mathbb{R}^{n \times n}} \;\; & \mdot{C}{X} \\
\mbox{subject to} \;\; & \mdot{A_i}{X} = b_i, \;\; i=1,\ldots m, \\
& X \succeq 0
\end{split}
\end{equation}
where $A_i \in \mathbb{R}^{n \times n}, b_i \in \mathbb{R}, i=1,\ldots,n$ and $C \in \mathbb{R}^{n \times n}$ are problem data, and $X$ is the optimization variable.
Most general purpose semidefinite programs typically use interior point methods \citep{JN06b} to solve problems to very high precision, but the methods  scale very poorly: interior point solvers operate in the space of $X\in\mathbb{R}^{n\times n}$, and the solution time is cubic in the number of variables, i.e., $O(n^6)$.

Owing to this fact, a number of approaches based upon \emph{low-rank factorization} have become popular for solving large-scale semi-definite programs.  Specifically, a now-classic result \citep{pataki1998rank} guarantees that for a problem with $m$ constraints,
there exists a rank-$\sqrt{2m}$ optimal solution.  Motivated by this fact, in a seminal work, \citet{burer2003nonlinear} proposed to solve the semidefinite program directly in terms of the low-rank factorization $X = V^T V$ for $V \in \mathbb{R}^{k \times n}$, where one can omit the semidefinite constraint (because it is implied by the factorization), and where $k$ could be chosen for instance to satisfy $k \geq \sqrt{2m}$; this leads to many fewer variables in the problem of interest. They specifically used an augmented Lagrangian method to solve the constrained optimization problem, and noticed that despite the fact that the resulting problem is of course nonconvex, in practice this would virtually always find the \emph{optimal} solution to the original SDP.  Although this remained an empirical rather than theoretical property for many years, very recently \citet{boumal2016non} has shown that there are \emph{no spurious local optima} in the augmented Lagrangian form for these problems, and thus such methods are guaranteed to find the optimal solution.

\subsection{The Mixing method}

Despite the appeal of the low-rank solution approach, the method is still relatively slow for even medium-sized SDPs; because augmented Lagrangian methods require a carefully-tuned gradient descent procedure to optimize the augmented Lagrangian and adjust the Lagrangian penalty parameter.  For a specific class of semidefinite program, however, \citet{wang2017mixing} developed a specific form of block coordinate descent, called the Mixing method, which substantially improves upon the generic low-rank approach both practically and theoretically.  Specifically, the Mixing method aims to solve the diagonally-constrained SDP
\begin{equation}
\minimize_{X \succeq 0 \in \mathbb{R}^{n \times n} } \; \mdot{C}{X}, \;\; \mbox{s.t.} \;\; X_{ii} = 1,\; i=1,\ldots,n
\end{equation}
i.e., an SDP where there are only $n$ equality constraints, constraining each diagonal entry of $X$ to be 1.  This is referred to as the MAXCUT SDP, since it corresponds to the SDP relaxation of the maximum cut problem, where $C$ is the adjacency matrix of a graph.  In factorized form $X = V^T V$, this corresponds to the non-convex problem
\begin{equation}
\minimize_{V \in \mathbb{R}^{k \times n}} \; \mdot{C}{V^T V}, \;\; \mbox{s.t.} \;\; \|v_i\|_2 = 1,\; i=1,\ldots,n,
\end{equation}
where $v_i \in \mathbb{R}^k$ denotes the $i$th column of $V$.  The approach is then quite simple: if we hold all but one $v_i$ fixed, there is a simple closed-form solution for the remaining $v_i$ given by
\begin{equation}
v_i := \mathrm{normalize}\left(-\sum_{j\neq i} c_{ij} v_j\right).
\end{equation}
The Mixing method simply repeats this iteration until convergence; unlike the augmented Lagrangian approach, the method does not require any step size or tuned penalty parameter.  In practice, it is also an order of magnitude or more faster than any other state-of-the-art solver for such problems.  And as shown by \citet{wang2017mixing}, as long as $k$ is chosen such that $k > \sqrt{2n}$, the method will still converge to the \emph{global} optimum of the original SDP.  This fast method for low-rank semidefinite programming will be the basis for our fast SDP-based MAX2SAT solver.

\section{MAXSAT and the semidefinite relaxation}

In this section we first derive a general SDP relaxation for the MAXSAT problem (although in the subsequent section we will focus on the MAX2SAT problem, we note that the formulation here is general and applies to \emph{any} MAXSAT problem, though it is substantially looser for settings with more than two variables per clause).  We then derive a version of the Mixing method specifically for this SDP; because the SDP is virtually identical to that for the MAXCUT problem, the algorithm is mathematically the same as the original Mixing method, with just a few additional elements introduced to make it computationally efficient for the specific problem structure of the MAXSAT.

\subsection{The MAXSAT SDP}
\begin{figure}[t!]\centering
\includegraphics[width=0.35\textwidth]{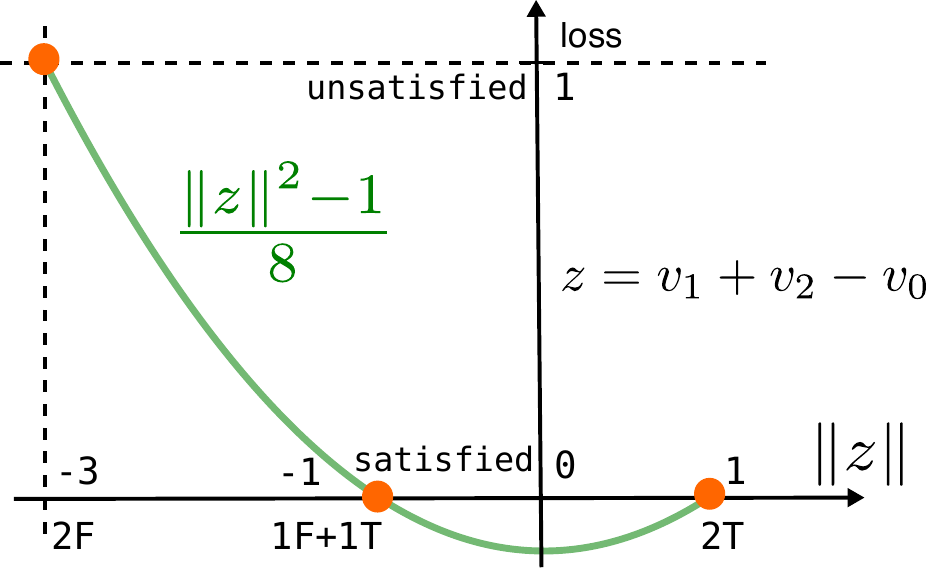}
\caption{Convex lower bound of the unsatisfiability loss.}
\label{fig:loss}
\end{figure}

For the purposes of this section, we will define the general MAXSAT problem slightly different than in the introduction, as the optimization problem 
\begin{equation}
  \maximize_{v \in\{-1,1\}^n}\;\; \sum_{j=1}^m \bigvee_{i=1}^n \mathbf{1}\{s_{ji} v_i > 0 \},
\end{equation}
where $v_i \in \{-1,1\}$ are the binary optimization variables (i.e., making these $+1/-1$ instead of 0/1 as we did in the introduction); and where $s_{ji} = \pm 1$ for variable in the clause (and taking each value depending on whether the variable is negated or not in that clause), and $s_{ji} = 0$ otherwise.  For conversion to SDP form it is slightly more convenient to formulate this in minimization, or \emph{unsatisfiability}, form
\begin{equation}
  \minimize_{V\in\{-1,1\}^n}\;\; \sum_{j=1}^m \bigwedge_{i=1}^n \mathbf{1}\{s_{ji} v_i < 0 \},
\end{equation}
where $\bigwedge$ is the logical and symbol.  At this point, we now seek a continuous \emph{lower-bound} on this objective, which we refer to simply as the \emph{loss}, and which takes the value $+1$ if an only if all variables within the clause are false (-1) and is zero or less otherwise, i.e.,
\begin{equation}
        \text{loss}(v; s_j) \leq \texttt{unsat}(v, s_j), \;\forall v_i\in\{-1,1\}^n.
\end{equation}
In order to translate the problem specifically to an SDP, we specifically search for a \emph{quadratic} loss; taking the \emph{greatest} such lower bound that satisfies the conditions above, we arrive at the loss function
\begin{equation}
        \text{loss}(v; s_j) = \frac{\left( \sum_{i=0}^n v_i s_{ji} \right)^2-(n_j-1)^2}{4n_j},
\end{equation}
where $n_j$ is the number of variables in clause $j$ (i.e, 2 in all cases for MAX2SAT), and where we define $v_0 = 1$ and $s_{j0} = -1$ (we introduce these variables to make the subsequent SDP relaxation easier to write, with purely quadratic terms in the objective).  This loss is best illustrated graphically, as shown in Figure \ref{fig:loss} for the case of $n_j = 2$.  For the case of any $n_j$, it's easy to verify that this quantity is equal to $+1$ if no clauses are satisfied, $0$ if exactly one or exactly $n_j$ clauses are satisfied, and is strictly less than 0 otherwise; these are exactly the conditions that the loss was required to meet.

We can now relax this binary problem into its vector form by relaxing the variable $v_i \in \{-1,+1\}$ to $v_i \in \mathbb{R}^k$ with $\|v_i\|_2 = 1$.  In this vector form, the loss above becomes simply
\begin{equation}
        \text{loss}(v; s_j) = \frac{\|V s_j\|^2 -(n_j-1)^2}{4n_j}
\end{equation}
which is equivalent to the SDP
\begin{equation}
\minimize_{V \in \mathbb{R}^{k \times n}} \; \mdot{S^T D^{-1} S}{V^T V}, \;\; \mbox{s.t.} \;\; \|v_i\|_2 = 1,\; i=1,\ldots,n,
\end{equation}
where the matrices are all as defined above and with $D_{jj} = 4n_j$.  This is precisely the form of the SDP solved by the Mixing method, with the only difference being the special form of the $C = S^T D^{-1} S$ matrix, which requires some slight optimizations to produce an efficient method, which we discuss next.  We also highlight that for the case of MAX2SAT, this is precisely the same form as the SDP relaxation of \citet{goemans1995improved}, but generalizes to longer clauses as well.

\subsection{The Mixing method for the MAXSAT SDP}
Recall from the previous section that the Mixing method for solving the SDP above reduces to the iterations
\begin{equation}
        v_i := \mathrm{normalize}\left(-\sum_{k\neq i} (S^T D^{-1} S)_{ik} v_k\right).
\end{equation}
To compute this term efficiently without requiring precomputation of the entire $S^T D^{-1} S$ matrix ahead of time (which is actually a \emph{dense} matrix owing to the fact that the ``0th'' variable $v_0$ participates in ever clause), note that we can cache the vectors $z_j \in \mathbb{R}^k$, $j=1,\ldots,m$, with
\begin{equation}
z_j = \sum_{i=0}^n s_{ji} v_i.
\end{equation}
Then the update equation can be written as
\begin{equation}
v_i := \mathrm{normalize}\left(-\sum_{j=1}^m \frac{s_{ji}}{4 n_j} z_j \right).
\end{equation}
if we first remove all the terms in each $z_j$ that depends on $v_i$.  The final form of the Mixing method for this MAXSAT SDP is given in Algorithm \ref{alg:mixing}.  Notice that updating the variable $v_i$ only requires time $O(k\cdot\texttt{nnz})$ times, where $\texttt{nnz}$ is the number of clauses that contain the $i$th variable.  Furthermore, although this would be a large number for the $v_0$ term, we can avoid updating $v_0$ entirely, as it plays the role of an arbitrary ``truth direction'' that need not be updated from its original setting.

\SetAlgoNoLine
\LinesNumbered
\begin{algorithm}[t]
        Initialize all $v_i$ randomly on a unit sphere\;
        Let $z_j=\sum_{i=0}^n s_{ji}v_i$ for $j=1,\ldots,m$\;
        \While{not yet converged}{
                \For{$i=1,\ldots,n$}{
                        \lForEach{$s_{ji}\neq 0$}{$z_j := z_j - s_{ji}v_i$}
                        $v_i := \textnormal{normalize}\left(-\sum_{j=1}^m \frac{s_{ji}}{4n_j} z_j\right)$\;
                        \lForEach{$s_{ji}\neq 0$}{$z_j := z_j + s_{ji}v_i$}
                }
        }
        \caption{The Mixing method for MAXSAT relaxation}\label{alg:mixing}
\end{algorithm}

\section{The MIXSAT algorithm}

In this section, we combine the previous SDP approach within a simple branch and bound framework to derive a exact and complete solver for MAX2SAT instances.\footnote{Since the aforementioned bound works for MAXSAT problems of any size, the methods here could be in theory applied to any MAXSAT problem. We will show that the methods generalize to MAX3SAT instances later in the paper, though theoretically the bound would be looser for longer clauses.}  Although the SDP relaxation and fast solution method plays the largest role in our solver, a number of other elements are crucial as well, such as warm starting based upon previous solutions in the tree, rounding to find feasible solutions at intermediate stages, and bounding the SDP optimal value by a particular set of primal and dual feasible points.  We also introduce data structures such as a watched stack (similar to those used by SAT solvers) to efficiently make inferences about those clauses that have been fulfilled or falsified.  Despite some additional complexity, we highlight that the ultimate algorithm here is still relatively simple compared to many state-of-the-art MAXSAT sovlers: just a branch and bound strategy with well-established optimizations adapted to this particular SDP-based heuristic.  The fact that the approach can nonetheless outperform existing solvers \emph{both} for complete and incomplete versions highlights the potential power of this SDP relaxation, not just from the theoretical side (which has been well-studied by past work), but from the practical side as well.

\SetAlgoNoLine
\LinesNumbered
\begin{algorithm}[t]
        Initialize a priority queue $Q = \{$initial problem$\}$\;
        \While{Q is not empty}{
                New SDP root $P$ = $Q$.pop()\;
                Solve $f^* := \texttt{SDP}(P)$ with the Mixing method\;
                \lIf{$\lceil f^*-\epsilon\rceil\geq \texttt{best}$}{continue}
                Update $\texttt{best}$ and resolution orders by randomized rounding on $V:=\arg \texttt{SDP}(P)$\;
                \ForEach{subproblems of $P$}{
                Initialize $\texttt{primal}$ and $\texttt{dual}$ objective values\;
                \lIf{$\texttt{primal}\leq \texttt{best}$}{Expand}
                \lElseIf{$\lceil\texttt{dual}\rceil \geq \texttt{best}$}{Prune}
                \lElse{Push the subproblem into the $Q$}
                }

        }
        \caption{The MIXSAT algorithm}\label{alg:mixsat}
\end{algorithm}

The basic method, which we refer to as MIXSAT (for Mixing applied to MAXSAT), is shown in Algorithm \ref{alg:mixsat}.  The basic approach proceeds like a generic branch and bound search, with a priority queue storing problem instances (by using a priority queue ordered by heuristic values versus a stack, we can recover either a best first or depth first search, which we use in the incomplete and complete versions of the solver respectively).  We pop these problems off the queue, solve the SDP via the Mixing method, and then round the result to potentially obtain a new candidate solution (the SDP solution also leads to a subsequent ordering of variables to split on).  For each split, we set the relevant set of variables be either true or false, and we initialize the primal and dual objectives of each new subproblem, as described below, then push these on to the queue and continue until the queue is empty (for the complete solver).  There are indeed several details of the process, which we describe next throughout this section.

\subsection{Pruning and rounding via SDP}

\paragraph{Pruning.}  Because our SDP approximation is a lower bound of the exact minimum unsatisfiability problem, it can serve as a lower-bound in the branch-and-bound framework.  Specifically, as is standard in branch and bound, if the objective value of a solved SDP value is higher than the current best known candidate solution, we can prune the entire subtree below that problem. Otherwise, we expand the subtree and push the subproblems on to the queue.

\paragraph{Rounding.}
The SDP solutions also serve as a good method for obtaining candidate solutions via rounding, which can usually very quickly provide a good best known value, even in early stages of the algorithm.   The randomized rounding method works by repetitively drawing random vectors $r$ uniformly from the sphere, and assigning binary $v_i=\text{sign}(v_0^T r\cdot v_i^T r)$ as a candidate binary solution. That is, if a variable $v_i$ is at the same side with the truth direction $v_0$, we assign it to be true ($+1$) and otherwise false ($-1$).  Finally, if the rounded solution ever obtains the same integer objective value as the best lower bound, we know we have found an optimal solution.

\subsection{Bounding SDP values by initializations}

Even with the pruning and rounding strategies above, it typically requires a large number of node expansions to solve even most MAX2SAT problem.  As one data point, for instance for the MAX2SAT problems with $120$ variables, it typically requires searching through about 10k nodes in order to provide a verified solution to the problem.  However, we reduce the number of expanded nodes (and hence the number of SDP solves) by reusing the solution of an SDP to quickly find both primal and dual feasible initialization to its subproblems.  Indeed, sometimes we don't even need to solve the SDP in order to expand or prune a subproblem.

\begin{figure}[t!]\centering
\includegraphics[width=0.25\textwidth]{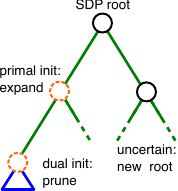}
\caption{Illustration of primal and dual initialization.}
\label{fig:tree}
\end{figure}

Recall the vector program equivalent to our diagonal constrained SDP relaxation
\begin{equation}
        \minimize_{V\in\bR^{k\times n}}\; f(V):=\mdot{C}{V^TV},\quad\text{s.t. }\;\norm{v_i}^2=1.
\end{equation}
By the SDP duality, we have the following dual problem
\begin{equation}
        \maximize_{\lambda\in\bR^{n}}\; D(\lambda):=-1^T\lambda,\quad\text{s.t. }\;C+D_{\lambda}\succeq 0.
\end{equation}
Any \emph{feasible} primal and dual solution will provide an upper and lower bound respectively of the optimal SDP objective value $f^*$.
Furthermore, we know that the best known discrete solution, $\texttt{BEST}$, is always larger than the optimal discrete solution $\texttt{UNSAT}$.
That is,
\begin{equation*}
        f(V) \geq f^* \geq D(\lambda) \quad\text{and}\quad \texttt{BEST}\geq \texttt{UNSAT} \geq \lceil f^*\rceil.
\end{equation*}
Thus, for any feasible primal solution $V$ and any feasible dual solution $\lambda$, we can prune or expand the subproblem when the following condition happens, as illustrated by Figure~\ref{fig:tree}.
\begin{itemize}
        \item $f(V)\leq \texttt{BEST}\Rightarrow$ Expand (not prunable by SDP)
        \item $\lceil D(\lambda)\rceil \geq \texttt{BEST}\Rightarrow$ Prune (not contain optimal solution)
\end{itemize}

\paragraph{Upper bound by primal initialization.}  The primal initialization is very simple: we simply copy all unassigned variables from the root solution into its subproblems, which will still provide a feasible solution.  The effective length of the remaining clauses do change when this occurs, which can be updated efficiently by maintaining each $z_j$ properly after such updates.

\paragraph{Lower bound by dual initialization.}
The dual initialization is more complex than the primal initialization because it involves the semidefinite $C+D_\lambda\succeq 0$ constraint.
The key here is that variable assignment can be viewed equivalently as moving the necessary coefficients from $s_{ji}$ to $s_{j0}$. For example,
assigning $v_1$ as false can be written as 
\begin{equation*}
        z_j = 1\cdot v_1+1\cdot v_2-1\cdot v_0 \;\;\xRightarrow{v_1:= F} \;\; z_j = 0\cdot v_1 + 1\cdot v_2 -2\cdot v_0.
\end{equation*}
Thus, if we look at the loss function, the coefficients $c_{ij}$ for unassigned $i$ and $j$ actually remains the same in the new cost matrix $\hat{C}$.
Thus, the difference in cost matrices $C$ and $\hat{C}$ after variable assignment will only reflects in the first column and row
\begin{equation}
        C-\hat{C} = \begin{pmatrix}0 & \delta^T \\ \delta & 0\end{pmatrix}
\end{equation}
If we add a proper vector $\xi=\begin{pmatrix}\norm{\delta}_1\\ |\delta|\end{pmatrix}$, we can see that
\begin{equation}
        C+D_{\lambda^*+\xi}=(C+D_{\lambda^*})+ \begin{pmatrix}\norm{\delta}_1& \delta^T \\ \delta & D_{|\delta|}\end{pmatrix}\succeq 0,
\end{equation}
where $\lambda^*$ is the optimal solution from the root.  That is, $\lambda^*+\xi$ is a feasible solution for the new dual problem, which provides our lower bound as desired.

We also note that given an \emph{optimal} solution of the SDP (as computed by Mixing), the optimal dual variable is simply given by
\begin{equation}
\label{eq:lambda}
\lambda_i = \|V c_i\|_2.
\end{equation}
Thus, when we solve the SDP optimally via the Mixing method, we also recover the optimal dual solution which serves as the starting point for future bounds.

\subsection{Data structures and implementation details}

Finally, there are a number of data structures and/or implementation details that are important for an efficient implementation of the branch and bound. For the full implementation and more technical details, please see the source code in \url{https://github.com/locuslab/mixsat} .

\paragraph{Variable ordering.}
We have found that the dual variable $\lambda$ in \eqref{eq:lambda} is also very useful in determining the resolution order (the order for splitting on variables to create the new subproblems).  
Thus, every time we solve an SDP, we will reorder the remaining variables according to the descending order of $\lambda$ (that is, variable with larger $\lambda_i$ will be resolved first).

\begin{figure}[t!]\centering
\includegraphics[width=0.3\textwidth]{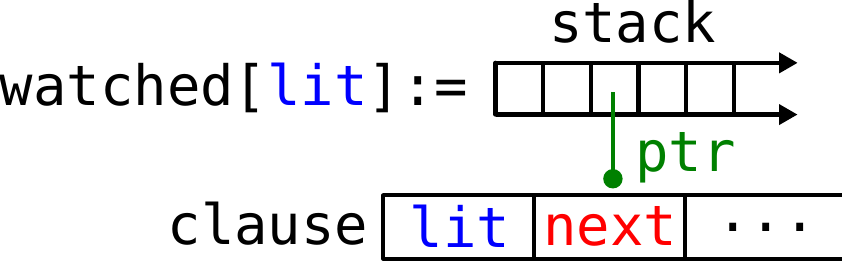}
\caption{Illustration the watched stack.}
\label{fig:watched}
\end{figure}

\paragraph{Inference via watched stack.}
In order to efficiently find the set of all clauses that contain a certain variable (including after variable reordering) and efficiently inferring which clauses have been fulfilled/falsified, we use a data structure similar to the watched literal technique in SAT solvers \citep{moskewicz2001chaff}.  Specifically, to accomplish this we used a \emph{watched stack} data structure.  For each literal, we maintain a stack of pointers, pointing to the next literal (in pivot order) in a clause.  The pointer is only pushed to the stack if it's the last pointer being able to falsify the clause. This way, assigning/freeing variables only touches those clauses that haven't been satisfied, and evaluating a set of variable assignment takes exactly $O(\texttt{nnz})$ time.
The watched stack is illustrated in Figure~\ref{fig:watched}.

\paragraph{Estimating SDP objective and balancing pruning/rounding.}
Because the Mixing method is iterative, it is usually desirable to only solve to some given precision, and not solve the SDPs exactly.
Since the Mixing method attains linear convergence \citep{wang2017mixing}, we estimate the converging constant by function decrease and estimate the distance to the optimal objective values.  Furthermore, we need to balance between the time spent in solving SDPs and the randomized roundings (because we rounds multiple times).  We have empirically found that taking rounding time proportion to the square root of the remaining variables works best. 

\paragraph{DFS vs BFS.}
We have described our branch-and-bound solver in an abstract manner in Algorithm~\ref{alg:mixsat}, where the priority queue could prioritize by any ordering.  When implementing the complete solver, where the full search tree needs to be searched anyway, we let the queue $Q$ be an stack; that is, we evaluate the nodes in the deep-first search (DFS) order.  For the incomplete version, we let $Q$ be the priority queue ordered by the clipped loss ($\sum_j \max(0, loss(V,s_j))$); that is, we explore the nodes in a best-first search (BFS) order.  However, the primal and dual initialization requires the DFS order to be efficiently implemented.  Thus, we run one depth-limited DFS for each of the BFS node evaluated even in the incomplete version.  Specifically, we limit the depth to be less than $8$ to limit the number of nodes (subproblems) pushed to the heap at each root.

\section{Experimental results}

\begin{figure}[t!]\centering
\includegraphics[width=0.5\textwidth]{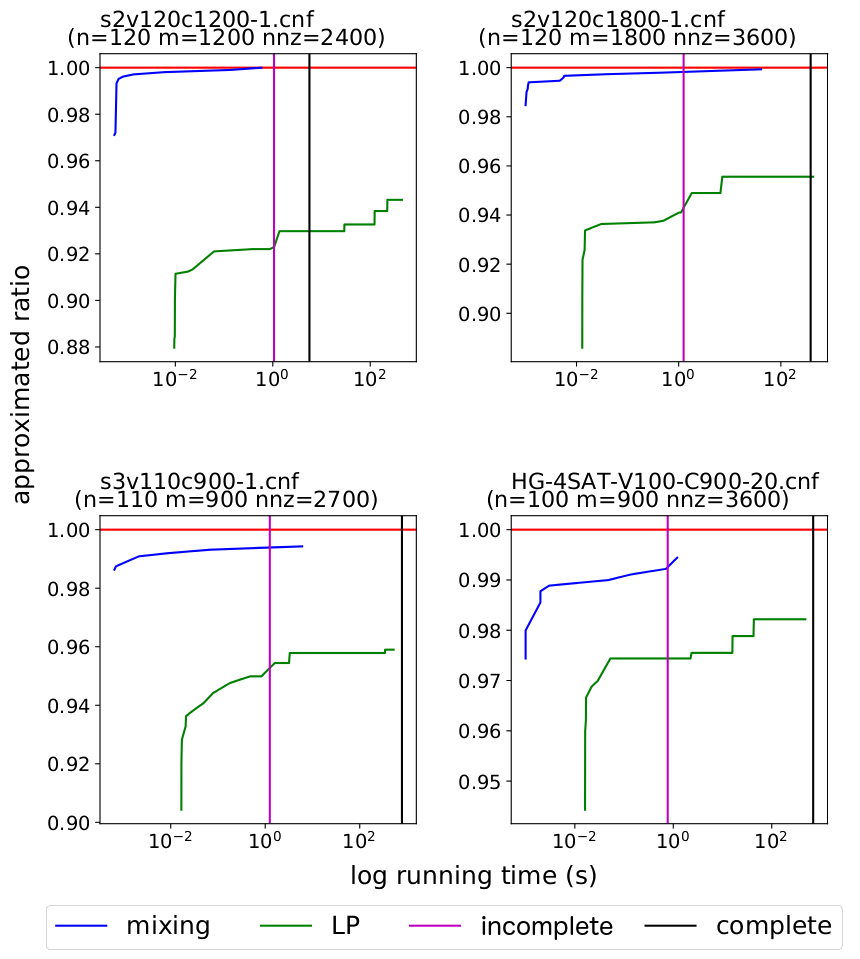}
\caption{The approximation ratio versus running time. We demonstrated that the Mixing method not only outputs solutions faster than LPs, but also admits higher approximation ratio.}
\label{fig:ratio}
\end{figure}
In this section we present the experimental results highlighting the performance of the MIXSAT method on the relevant problems from the MAXSAT 2016 competition (later years have eliminated the ``random'' track of the competition, and thus don't have a sufficient number of MAX2SAT problems).  Specifically, we test our solvers on all 228 the MAX2SAT instances (s2v instances) from the unweighted random categories.  As a comparison point, for each problem we selected the \emph{best} performing entry for both the complete and incomplete tracks in the MAXSAT 2016 contest; since different solvers often performed best on different problems, this sets a very high bar for comparison, though with the large caveat mentioned earlier that of course the other solvers were attempting to optimize performance over all tasks, not just MAX2SAT problems.  Nonetheless, given the level of specializing in existing approaches, we feel that the comparison provides strong evidence that our SDP can achieve compelling performance in this domain. 
Further, to test the performance on harder problems, we evaluated our solver on the 107 crafted MAXCUT instances in the 2016/2018 competitions, comparing it to the best solvers in 2016/2018.
We also tested our solver on all the random MAX3SAT instances in 2016 to see how it generalized to different clause length.
For fair comparison, we replicate the experimental setup of the 2016 environment, with exactly the same CPU (Intel Xeon E5-2620) in single core mode and a 3.5 GB memory limit. 

\paragraph{The approximation ratio of MAXSAT SDPs.}
As a starting point, we first highlight the power of the MAXSAT Mixing solver combined just with randomized rounding (i.e., just looking at a single relaxation, without any branch and bound).  These results are shown in Figure \ref{fig:ratio}.  The starting points of the blue and green curves indicate the time that the Mixing / LP solver (by Gurobi) output the solution, and the following curves indicate the highest approximation ratio generated by the LP/SDP randomized rounding procedures.  We can see that the Mixing method actually output solution faster than the LP solver,  and has much higher approximation ratio than LPs.  Indeed, sometimes the simple randomized rounding procedure can output the optimal solution faster than the best complete and incomplete solvers (vertical black and purple lines).

\begin{figure}[t!]\centering
\includegraphics[width=0.5\textwidth]{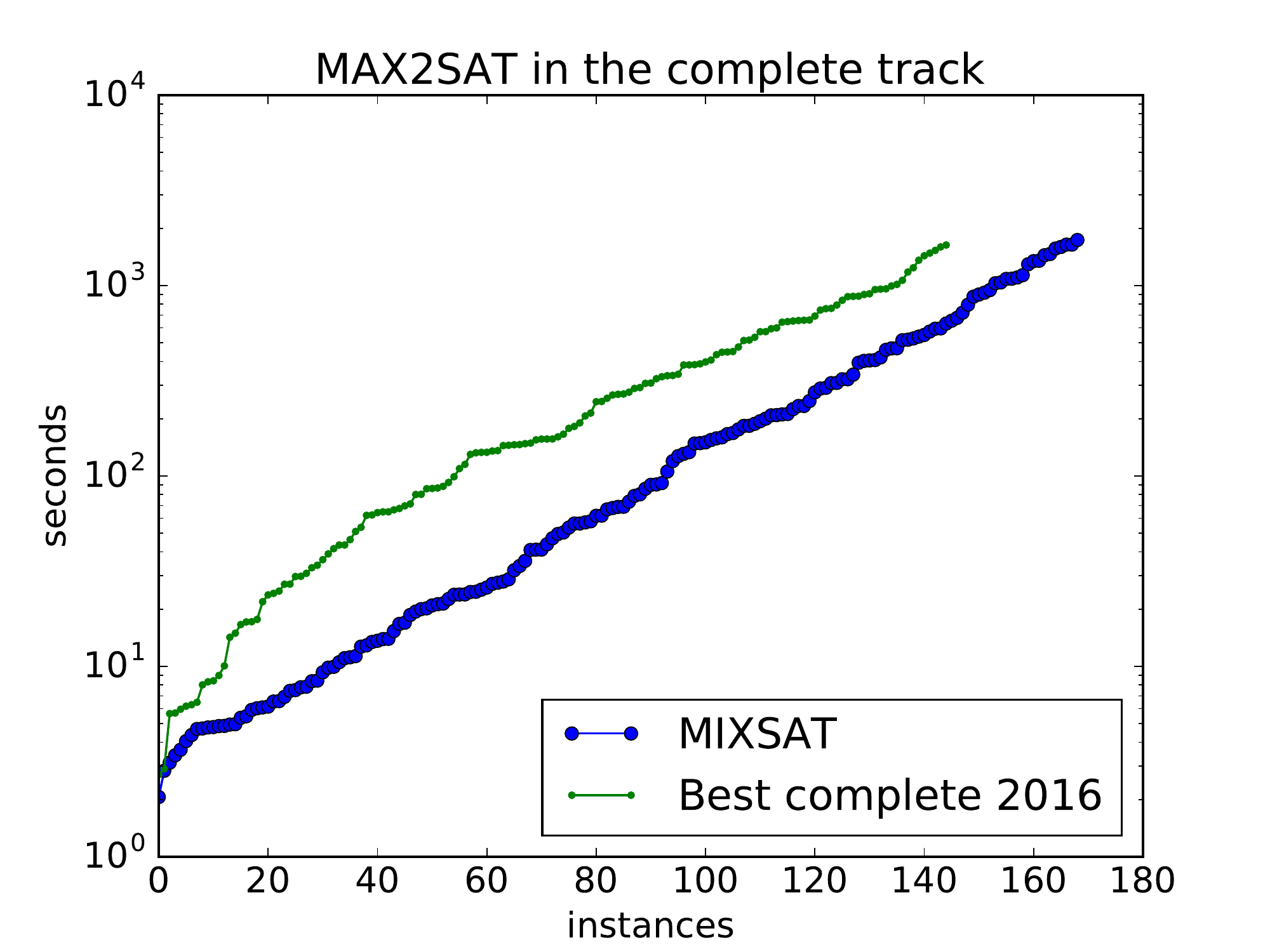}
\caption{Experimental result of MAX2SAT in the complete track. We demonstrate that the MIXSAT algorithm consistently outperforms the best complete solvers in the MAXSAT 2016 competition in every MAX2SAT instances.
        Specifically, we solved 169/228 instances in 30 mins in avg 273 secs, while the best solvers solved only 145/228 instances in avg 341 secs.}
\label{fig:complete}
\end{figure}

\begin{figure}[t!]\centering
\includegraphics[width=0.5\textwidth]{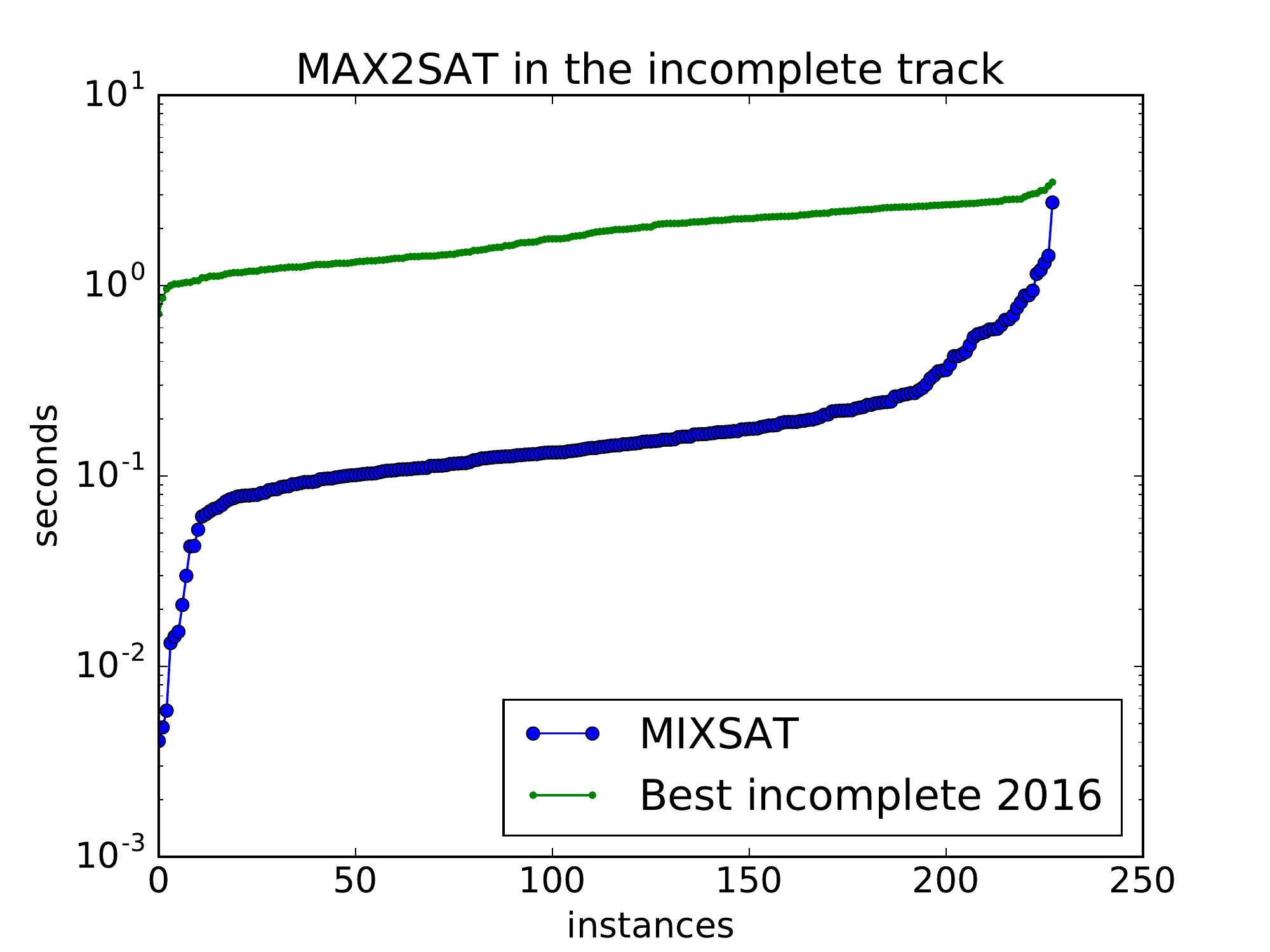}
\caption{Experimental result of MAX2SAT in the incomplete track. We demonstrate that the MIXSAT algorithm outperform the best incomplete solvers in the MAXSAT 2016 competition in every MAX2SAT instances.
        The MIXSAT algorithm used 0.22 sec in average to obtain the optimal solution, while the best solvers use 1.92 sec in average. That is, we are approximately 8.72 times faster than the best solvers.}
\label{fig:incomplete}
\end{figure}

\paragraph{Complete solvers.}
Next, we evaluate our solver on the above mentioned 228 MAX2SAT instances using the complete track rules, which only allows to output verified optimal solution and has a 30 minutes time limit.  The result is presented in Figure~\ref{fig:complete}, showing the number of solved problems versus time for MIXSAT shown against the best MAXSAT 2016 solver results.  We solved \textbf{169} out of the 228 MAX2SAT instances in the random categories for the complete solvers in an average of 273 sec, while the best solvers solved only \textbf{145} instances in an average of 341 seconds (note however that the average times are not truly comparable here as we solve substantially more problems).

\paragraph{Incomplete solvers.}
We test our solvers on the 228 MAX2SAT instances using the incomplete track rules, which allows the solver to output solutions before verification and has a 5 minutes time limit. The result is presented in Figure~\ref{fig:incomplete}. We recovered the optimal solution of all the 228 MAX2SAT instances in the random categories for the incomplete solvers in an average of 0.22 sec, while the best solvers solved in averagely 1.92 sec. That is, we are \textbf{8.72x faster} than the best solvers in average. To be specific, we found that the 2016 best solvers we compared with are usually one of CCLS \citep{luo2015ccls}, CnC-LS \cite{luo2017ccehc}, SC2016 \citep{wagnermaxsat}, and Swcca-ms \citep{cai2013local}. 

\begin{figure}[t!]\centering
\includegraphics[width=0.5\textwidth]{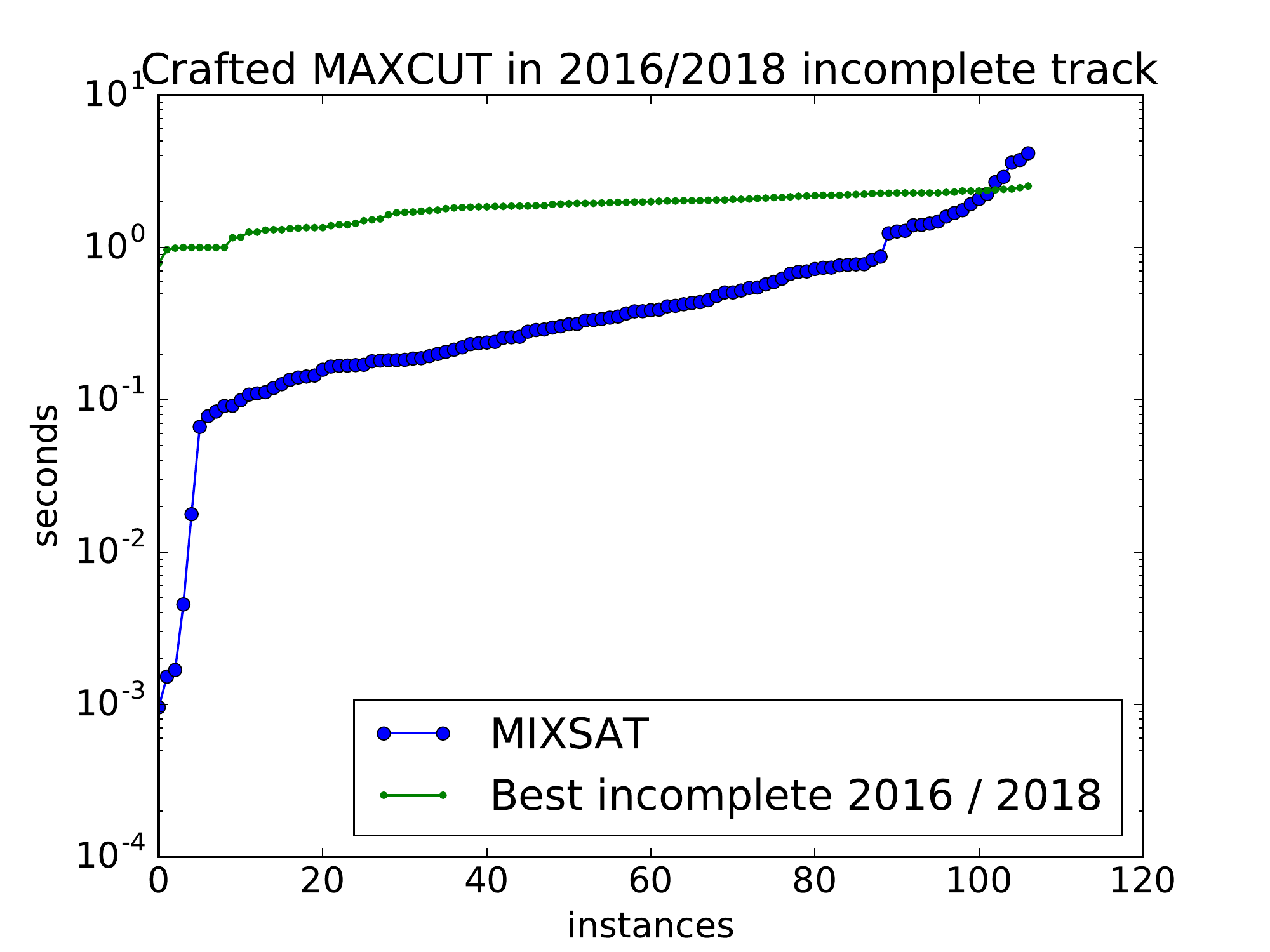}
\caption{Experimental result of the crafted MAXCUT instances in the 2016/2018 incomplete track. 
        The MIXSAT algorithm used 0.62 sec in average to obtain the optimal solution, while the best solvers in 2016/2018 took 1.84 sec in average. That is, we are still 3x faster in the crafted tracks.}
\label{fig:crafted}
\end{figure}

\begin{figure}[t!]\centering
\includegraphics[width=0.5\textwidth]{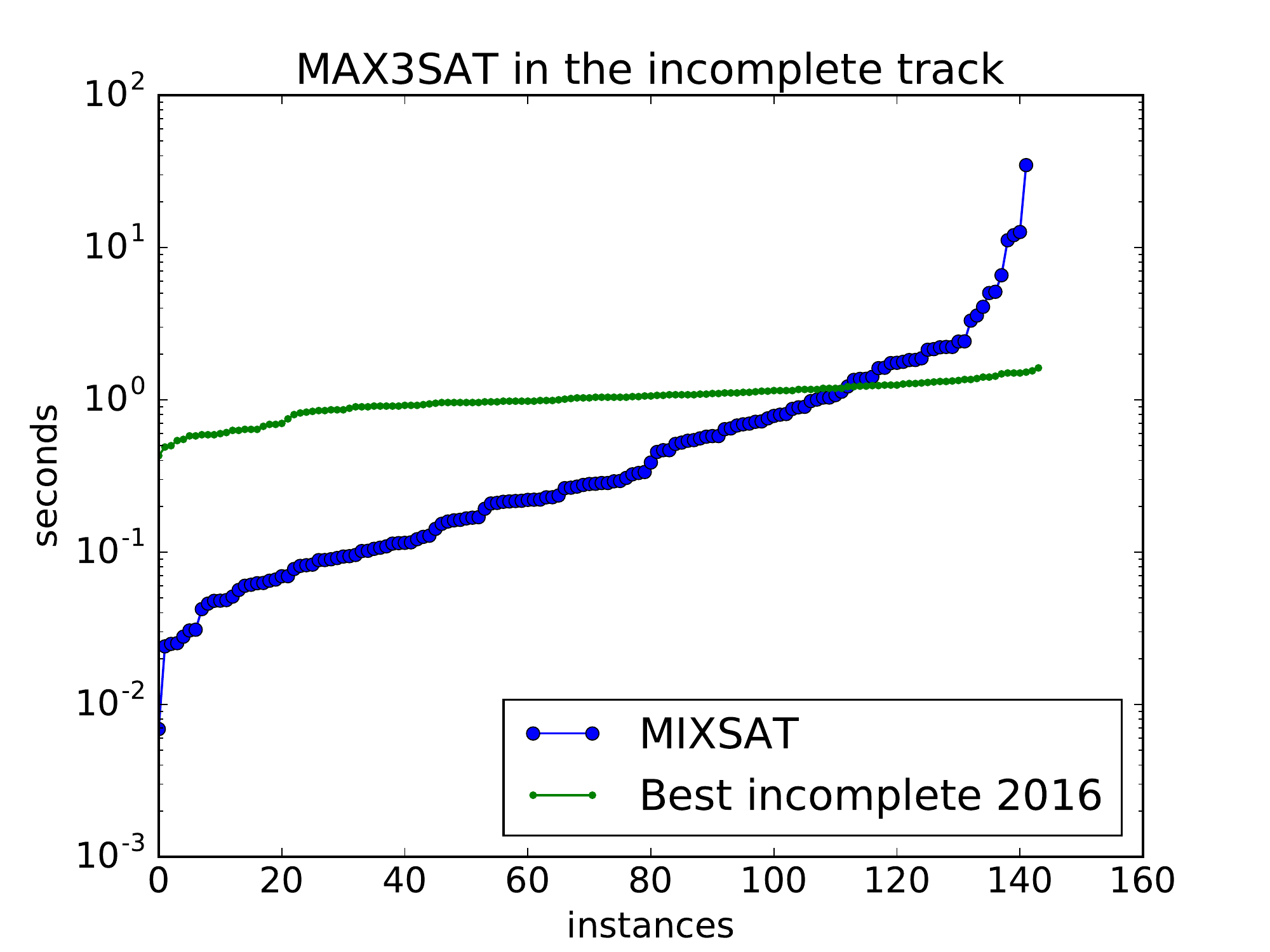}
\caption{Experimental result of MAX3SAT in the incomplete random track. Even though our solver is not designed for MAX3SAT,
       it still outperformed the 2016 best solvers in 114/144 instances and recover all the optimal solutions except for 2 instances.}
\label{fig:max3sat}
\end{figure}

\paragraph{Crafted MAXCUT and random MAX3SAT instances.}
Finally, we evaluate our solver against the 107 crafted MAXCUT instances in the 2016/2018 completions, as shown in Figure~\ref{fig:crafted}.
The result demonstrates our solver also performed well in harder instances: we recovered all optimal solution in avg 0.63 sec, which is \textbf{3x faster} comparing to the best solvers in 2016/2018 (avg 1.92s). 
Out of curiosity, we also tested our solver on the 144 random MAX3SAT instances in the 2016 competition. The result, shown in Figure~\ref{fig:max3sat}, suggests that our solver is still competitive even if it is not designed for MAX3SAT. It outperformed the best solvers in 2016 on 114/144 of the instances, though didn't recover the optimal solution for 2 instances.

\section{Conclusion}

This paper has highlighted that an approach, based upon low-rank semidefinite programming and a simple branch-and-bound strategy, is extremely competitive for solving the MAX2SAT problem, compared with some state-of-the-art solvers in previous contest on MAXSAT solving. Further, preliminary experiments demonstrate the performance also generalizes to the MAX3SAT problem.  Although the MAX2SAT/MAX3SAT domain is obviously a limited setting, the fact that an approach based upon SDP relaxations can perform well here offers substantial evidence that with modern methods for SDP solving, such relaxations can finally go beyond merely theoretically interesting, but become practical tools in the toolset of combinatorial optimization. These results suggest that they could be a powerful tool (one of many) to integrate into current solvers.  More fundamentally, we hope that the work serves as a starting point for more continuous approaches integrated into combinatorial search, ultimately providing a powerful new avenue for some fundamental problems in AI.

\bibliographystyle{aaai}
\bibliography{sdp.bib}

\appendix

\end{document}